\documentclass[reqno,12pt]{amsart}
\usepackage{amssymb}
\usepackage{amsfonts}
\usepackage{amsmath}

\newtheorem{theorem}{Theorem}
\theoremstyle{plain}
\newtheorem{acknowledgement}{Acknowledgement}

\newtheorem{corollary}{Corollary}

\newtheorem{definition}{Definition}

\newtheorem{remark}{Remark}

\numberwithin{equation}{section}

\input{tcilatex}

\begin{document}
\author{}
\title{}
\maketitle

\begin{center}
\thispagestyle{empty} \pagestyle{myheadings} 
\markboth{\bf  Yilmaz Simsek,
}{\bf Analysis of the p-adic q-Volkenborn Integrals... }

\textbf{\Large Analysis of the }$\mathbf{\Large p}$\textbf{\Large -adic }$%
\mathbf{\Large q}$\textbf{\Large -Volkenborn Integrals: an approach to
Apostol-type special numbers and polynomials}

\bigskip

\textbf{Yilmaz Simsek}

\bigskip

\textit{Department of Mathematics, Faculty of Science University of Akdeniz
TR-07058 Antalya, Turkey, }

\textit{E-mail: ysimsek@akdeniz.edu.tr\\[0pt]
}

\bigskip

\textbf{{\large {Abstract}}}\medskip
\end{center}

By applying the $p$-adic $q$-Volkenborn Integrals including the bosonic and
the fermionic $p$-adic integrals on $p$-adic integers,\textbf{\Large \ }we
define generating functions, attached to the Dirichlet character, for the\
generalized Apostol-Bernoulli numbers and polynomials, the\ generalized
Apostol-Euler numbers and polynomials, generalized Apostol-Daehee numbers
and polynomials, and also\ generalized Apostol-Changhee numbers and
polynomials. We investigate some properties of these numbers and polynomials
with their generating functions. By using these generating functions and
their functional equation, we give some identities and relations including
the generalized Apostol-Daehee and Apostol-Changhee numbers and polynomials,
the Stirling numbers, the Bernoulli numbers of the second kind,
Frobenious-Euler polynomials, the generalized Bernoulli numbers and the
generalized Euler numbers and the Frobenious-Euler polynomials. By using the
bosonic and the fermionic $p$-adic integrals, we derive integral
represantations for the generalized Apostol-type Daehee numbers and the\
generalized Apostol-type Changhee numbers.

\begin{quotation}
\bigskip
\end{quotation}

\noindent \textbf{2010 Mathematics Subject Classification.} 11B68; 05A15;
05A19; 12D10; 26C05; 30C15.

\noindent \textbf{Key Words.} the $p$-adic $q$-Volkenborn Integrals,
Dirichlet character, generalized Apostol-Bernoulli numbers and polynomials,
the\ generalized Apostol-Euler numbers and polynomials, generalized
Apostol-Daehee numbers and polynomials,\ generalized Apostol-Changhee
numbers and polynomials, Stirling numbers, Frobenious-Euler polynomials,
Generating functions.

\bigskip

\section{Introduction}

In recent years, special numbers and polynomials have played role in almost
all areas of mathematics, in mathematical physics, in computer science, in
engineering problems and in other areas of science. There are many different
methods and techniques for investigating and examining the special numbers
and polynomials including the Bernoulli numbers and polynomials, the Euler
numbers and polynomials, the Stirling numbers and the Apostol-type numbers
and polynomials, the Frobenious-Euler polynomials and others. The main
motivation of this paper is to construct generating functions for new
families of numbers and polynomials including genearlized Apostol-Daehee
numbers and genearlized Apostol-Changhee numbers attached to Dirichlet
character $\chi $. For constructing these functions, we can use the $p$-adic 
$q$-Volkenborn integrals techniques including integral equations of the
bosonic and the fermionic $p$-adic integrals on $p$-adic integers. These
numbers and polynomials are related to well-known numbers and polynomials
such as the generalized Bernoulli and Euler numbers and polynomials, the
Apostol-Bernoulli numbers and polynomials and the Apostol-Euler numbers and
polynomials, the Stirling numbers of the first and of the second kind, the
Daehee numbers and polynomials and also the Changhee numbers and
polynomials. By using the integral equations of the bosonic and the
fermionic $p$-adic integrals and generating functions, we give some integal
formulas and identities including the Lah numbers, the falling and rising
factorials, the well-known special numbers and polynomials and combinatorial
sums.

Throughout of this paper we need the following notations, definitions and
some families of the special numbers and polynomials.

The Apostol-Bernoulli polynomials, $\mathcal{B}_{n}(x;\lambda )$ are defined
by means of the following generating function: 
\begin{equation}
F_{A}(t,x;\lambda )=\frac{t}{\lambda e^{t}-1}e^{tx}=\sum_{n=0}^{\infty }%
\mathcal{B}_{n}(x;\lambda )\frac{t^{n}}{n!},  \label{Ap.B}
\end{equation}%
where $\lambda $ is a complex number and ($\left\vert t\right\vert <2\pi $
when $\lambda =1$ and $\left\vert t\right\vert <\left\vert \log \lambda
\right\vert $ when $\lambda \neq 1$) with%
\begin{equation*}
\mathcal{B}_{n}(\lambda )=\mathcal{B}_{n}(0;\lambda )
\end{equation*}%
denotes so-called $\lambda $-Bernoulli numbers (\textit{cf}. \cite%
{Kim2006TMIC}, \cite{KIMjang}, \cite{Luo}, \cite{ozden}, \cite%
{Srivastava2011}, \cite{srivas18}; see also the references cited in each of
these earlier works). Observe that%
\begin{equation*}
B_{n}=\mathcal{B}_{n}(0;1)
\end{equation*}%
denotes the Bernoulli polynomials of the first kind\textit{\ }(\textit{cf}. 
\cite{Bayad}-\cite{Wang}; see also the references cited in each of these
earlier works).

Kim \textit{et al}. \cite{TkimJKMS} defined the $\lambda $-Bernoulli
polynomials (Apostol-type Bernoulli polynomials), $\mathfrak{B}%
_{n}(x;\lambda )$ by means of the following generating function:%
\begin{equation}
F_{B}(t,x;\lambda )=\frac{\log \lambda +t}{\lambda e^{t}-1}%
e^{tx}=\sum_{n=0}^{\infty }\mathfrak{B}_{n}(x;\lambda )\frac{t^{n}}{n!}
\label{laBN}
\end{equation}%
($\left\vert t\right\vert <2\pi $ when $\lambda =1$ and $\left\vert
t\right\vert <\left\vert \log \lambda \right\vert $ when $\lambda \neq 1$)
with%
\begin{equation*}
\mathfrak{B}_{n}(\lambda )=\mathfrak{B}_{n}(0;\lambda )
\end{equation*}%
denotes the $\lambda $-Bernoulli numbers (Apostol-type Bernoulli numbers) (%
\textit{cf}. \cite{TkimJKMS}, \cite{jandY1}, \cite{srivas18}, \cite%
{Simsek2019}, \cite{simsek2017ascm}). A few of the $\lambda $-Bernoulli
numbers are given as follows:%
\begin{equation*}
\mathfrak{B}_{0}(\lambda )=\frac{\log \lambda }{\lambda -1},
\end{equation*}%
and%
\begin{equation*}
\mathfrak{B}_{1}(\lambda )=\frac{\lambda -1-\lambda \log \lambda }{\left(
\lambda -1\right) ^{2}}.
\end{equation*}

If $n>1$, than we have%
\begin{equation*}
\mathfrak{B}_{n}(\lambda )=\lambda \left( \mathfrak{B}_{n}(\lambda
)+1\right) ^{n}.
\end{equation*}

Kim \textit{et al}. gave a relation between these numbers and the
Frobenius-Euler numbers \cite[Theorem 1, p. 439]{TkimJKMS} as follows:%
\begin{equation*}
\mathfrak{B}_{0}(\lambda )=\frac{\log \lambda }{\lambda -1}H_{0}\left( \frac{%
1}{\lambda }\right)
\end{equation*}%
and%
\begin{equation*}
\mathfrak{B}_{n}(\lambda )=\frac{\log \lambda }{\lambda -1}H_{n}\left( \frac{%
1}{\lambda }\right) +\frac{n}{\lambda -1}H_{n-1}\left( \frac{1}{\lambda }%
\right)
\end{equation*}%
where $H_{n}\left( \frac{1}{\lambda }\right) $ denotes the Frobenius-Euler
numbers, defined by means of the following generating function:

Let $u$ be a complex numbers with $u\neq 1$.%
\begin{equation}
F_{f}(t,u)=\frac{1-u}{e^{t}-u}=\sum_{n=0}^{\infty }H_{n}(u)\frac{t^{n}}{n!}
\label{FEN}
\end{equation}%
(\textit{cf}. \cite{DSKIMfrob}, \cite[Theorem 1, p. 439]{TkimJKMS}, \cite%
{jnt2015}, \cite{srivas18}; see also the references cited in each of these
earlier works).

Let $r$ be a positive integer, and let $\lambda \neq 1$ be any nontrivial $r$%
-th root of $1$. Observe that the $\lambda $-Bernoulli numbers are reduced
to the Apostol-Bernoulli numbers and the twisted Bernoulli numbers ( \textit{%
cf}. \cite[Theorem 1, p. 439]{TkimJKMS}; see also the references cited in
each of these earlier works).

We \cite{ysimsek Ascm} gave the following functional equation:%
\begin{equation*}
F_{B}(t,0;\lambda )=\frac{\log \lambda }{\lambda -1}F_{f}\left( t,\frac{1}{%
\lambda }\right) +F_{A}(t,0;\lambda ).
\end{equation*}%
By using the above functional equation, we have%
\begin{equation*}
\mathfrak{B}_{n}(\lambda )=\frac{\log \lambda }{\lambda -1}H_{n}\left( \frac{%
1}{\lambda }\right) +\mathcal{B}_{n}(\lambda )
\end{equation*}%
(\textit{cf}. \cite{ysimsek Ascm}).

Let $\chi $ be a non-trivial Dirichlet character with conductor $d$. Let $%
\lambda $ be a complex number. The generalized Apostol-Bernoulli numbers
attached to Dirichlet character,\ $\mathcal{B}_{n,\chi }(\lambda )$ are
defined by means of the following generating function:%
\begin{equation}
\sum_{j=0}^{d-1}\frac{\lambda ^{j}e^{tj}t\chi (j)}{\lambda ^{d}e^{td}-1}%
=\sum_{n=0}^{\infty }\mathcal{B}_{n,\chi }(\lambda )\frac{t^{n}}{n!},
\label{GA}
\end{equation}%
where $\left\vert t+\log \lambda \right\vert <\frac{2\pi }{d}$ (\textit{cf}. 
\cite{apostol}, \cite{MSKIM}, \cite{KIMaml2008}, \cite{KIMjmaa2017}, \cite%
{RyooCHARbernoul}, \cite{srivas18}; see also the references cited in each of
these earlier works).

By combining (\ref{GA}) with (\ref{Ap.B}), we have%
\begin{equation*}
\mathcal{B}_{n,\chi }(\lambda )=d^{n-1}\sum_{j=0}^{d-1}\lambda ^{j}\chi (j)%
\mathcal{B}_{n}\left( \frac{j}{d};\lambda ^{p}\right) 
\end{equation*}%
for the trivial character $\chi \equiv 1$, we have%
\begin{equation*}
\mathcal{B}_{n}(\lambda )=\mathcal{B}_{n,1}(\lambda )
\end{equation*}%
(\textit{cf}. \cite{apostol}, \cite{MSKIM}, \cite{KIMaml2008}, \cite%
{KIMjmaa2017}, \cite{RyooCHARbernoul}, \cite{srivas18}).

The Apostol-Euler polynomials of first kind, $\mathcal{E}_{n}(x,\lambda )$
are defined by means of the following generating function:%
\begin{equation}
F_{P1}(t,x;k,\lambda )=\frac{2}{\lambda e^{t}+1}e^{tx}=\sum_{n=0}^{\infty }%
\mathcal{E}_{n}(x,\lambda )\frac{t^{n}}{n!},  \label{Cad3}
\end{equation}%
($\left\vert t\right\vert <\pi $ when $\lambda =1$ and $\left\vert
t\right\vert <\left\vert \ln \left( -\lambda \right) \right\vert $ when $%
\lambda \neq 1$), $\lambda \in \mathbb{C}$. Substituting $x=0$ into (\ref%
{Cad3}), we have The first kind Apostol-Euler numbers of order $k$:%
\begin{equation*}
\mathcal{E}_{n}(\lambda )=\mathcal{E}_{n}(0,\lambda )
\end{equation*}

Setting $\lambda =1$ into (\ref{Cad3}), one has the first kind Euler numbers%
\begin{equation*}
E_{n}=\mathcal{E}_{n}^{(1)}(1)
\end{equation*}%
(\textit{cf}. \cite{Bayad}-\cite{SrivastavaLiu}; see also the references
cited in each of these earlier works).

Let $\chi $ be a non-trivial Dirichlet character with conductor $d$. Let $%
\lambda $ be a complex number. The generalized Apostol-Euler numbers
attached to Dirichlet character,\ $\mathcal{E}_{n,\chi }(x,\lambda )$ are
defined by means of the following generating function:

\begin{equation}
2\sum_{j=0}^{d-1}\frac{\lambda ^{j}e^{tj}\chi (j)}{\lambda ^{d}e^{td}+1}%
=\sum_{n=0}^{\infty }\mathcal{E}_{n,\chi }(\lambda )\frac{t^{n}}{n!}
\label{GE}
\end{equation}%
where $\left\vert t+\log \lambda \right\vert <\frac{\pi }{d}$ (\textit{cf}. 
\cite{KIMaml2008}, \cite{KIMjmaa2017}, \cite{srivas18}; see also the
references cited in each of these earlier works).

By combining (\ref{GE}) with (\ref{Cad3}), we have%
\begin{equation*}
\mathcal{E}_{n,\chi }(\lambda )=d^{n}\sum_{j=0}^{d-1}\lambda ^{j}\chi (j)%
\mathcal{E}_{n}\left( \frac{j}{d};\lambda ^{p}\right)
\end{equation*}%
for the trivial character $\chi \equiv 1$, we have%
\begin{equation*}
\mathcal{E}_{n}(\lambda )=\mathcal{E}_{n,1}(\lambda )
\end{equation*}%
(\textit{cf}. \cite{KIMaml2008}, \cite{KIMjmaa2017}, \cite{srivas18}).

The Stirling numbers of the first kind, $S_{1}(n,k)$ are defined by means of
the following generating function:%
\begin{equation}
F_{S1}(t,k)=\frac{\left( \log (1+t)\right) ^{k}}{k!}=\sum_{n=0}^{\infty
}S_{1}(n,k)\frac{t^{n}}{n!}.  \label{S1}
\end{equation}%
These numbers have the following properties:

$S_{1}(0,0)=1$. $S_{1}(0,k)=0$ if $k>0$. $S_{1}(n,0)=0$ if $n>0$. $%
S_{1}(n,k)=0$ if $k>n$ and also%
\begin{equation*}
S_{1}(n+1,k)=-nS_{1}(n,k)+S_{1}(n,k-1)
\end{equation*}%
(\textit{cf}. \cite{Charamb}, \cite{Bayad}, \cite{Chan}, \cite{Roman}, \cite%
{SimsekFPTA}, \cite{AM2014}; and see also the references cited in each of
these earlier works).

The Bernoulli polynomials of the second kind, $b_{n}(x)$ are defined by
means of the following generating function:%
\begin{equation}
F_{b2}(t,x)=\frac{t}{\log (1+t)}(1+t)^{x}=\sum_{n=0}^{\infty }b_{n}(x)\frac{%
t^{n}}{n!}  \label{Br2}
\end{equation}%
(\textit{cf}. \cite[PP. 113-117]{Roman}; see also the references cited in
each of these earlier works).

The Bernoulli numbers of the second kind $b_{n}(0)$ are defined by means of
the following generating function: 
\begin{equation}
F_{b2}(t)=\frac{t}{\log (1+t)}=\sum_{n=0}^{\infty }b_{n}(0)\frac{t^{n}}{n!}.
\label{Be-1t}
\end{equation}%
These numbers are computed by the following formula:%
\begin{equation*}
\sum_{k=0}^{n-1}(-1)^{k}\left( 
\begin{array}{c}
n \\ 
k%
\end{array}%
\right) b_{k}(0)=n!\delta _{n,1},
\end{equation*}%
where $\delta _{n,1}$ denotes\ the Kronecker delta (\textit{cf}. \cite[p. 116%
]{Roman}). The Bernoulli polynomials of the second kind are defined by%
\begin{equation*}
b_{n}(x)=\int_{x}^{x+1}(u)_{n}du.
\end{equation*}%
Substituting $x=0$ into the above equation, one has%
\begin{equation*}
b_{n}(0)=\int_{0}^{1}(u)_{n}du.
\end{equation*}%
The numbers $b_{n}(0)$ are also so-called the \textit{Cauchy numbers} (%
\textit{cf}. \cite[p. 116]{Roman}, \cite{komotsu}, \cite{TKimTAKAO}, \cite%
{Qi}, \cite{ysimsek Ascm}, \cite{Roman}; see also the references cited in
each of these earlier works). In \cite{Kim2016Bernoull2}, Kim \textit{et al.}
gave a computation method for the Bernoulli polynomials of the second kind
are defined as follows:%
\begin{equation*}
b_{n}(x)=\sum_{l=0}^{n}\frac{S_{1}(n,l)}{l+1}\left(
(x+1)^{l+1}-x^{l+1}\right)
\end{equation*}%
and also Roman \cite[p.115]{Roman} gave%
\begin{equation*}
b_{n}(x)=b_{n}(0)+\sum_{l=1}^{n}\frac{nS_{1}(n-l,l-1)}{l}x^{l}.
\end{equation*}

By using the above formula for the Bernoulli polynomials and numbers of the
second kind, few of these numbers are computed as follows, respectively:

\begin{eqnarray*}
b_{0}(x) &=&1, \\
b_{1}(x) &=&\frac{1}{2}(2x+1), \\
b_{2}(x) &=&\frac{1}{6}(6x^{2}-1), \\
b_{3}(x) &=&\frac{1}{4}(4x^{3}-6x^{2}1), \\
b_{4}(x) &=&\frac{1}{30}(30x^{4}-120x^{3}+120x^{2}-19),\ldots
\end{eqnarray*}

and%
\begin{equation*}
b_{0}(0)=1,b_{1}(0)=\frac{1}{2},b_{2}(0)=-\frac{1}{12},b_{3}(0)=\frac{1}{24}%
,b_{4}(0)=-\frac{19}{720},\ldots .
\end{equation*}

In \cite{SimsekFPTA}, Simsek defined the $\lambda $-array polynomials $%
S_{v}^{n}(x;\lambda )$ by the following generating function:%
\begin{equation}
F_{A}(t,x,v;\lambda )=\frac{\left( \lambda e^{t}-1\right) ^{v}}{v!}%
e^{tx}=\sum_{n=0}^{\infty }S_{v}^{n}(x;\lambda )\frac{t^{n}}{n!},
\label{ARY-1}
\end{equation}%
where $v\in \mathbb{N}_{0}$ and $\lambda \in \mathbb{C}$ (\textit{cf}. \cite%
{Bayad}, \cite{Chan}, \cite{SimsekFPTA}, \cite{AM2014}; and see also the
references cited in each of these earlier works).

The $\lambda $-Stirling numbers, $S_{2}(n,v;\lambda )$ are defined by means
of the following generating function:%
\begin{equation}
F_{S}(t,v;\lambda )=\frac{\left( \lambda e^{t}-1\right) ^{v}}{v!}%
=\sum_{n=0}^{\infty }S_{2}(n,v;\lambda )\frac{t^{n}}{n!},  \label{SN-1}
\end{equation}%
where $v\in \mathbb{N}_{0}$ and $\lambda \in \mathbb{C}$ (\textit{cf}. \cite%
{Luo}, \cite{SimsekFPTA}, \cite{Srivastava2011}; see also the references
cited in each of these earlier works). By using (\ref{SN-1}), one easily
compute the following values for $S_{2}(n,v;\lambda )$:%
\begin{equation*}
S_{2}(0,0;\lambda )=1,S_{2}(1,0;\lambda )=0,S_{2}(1,1;\lambda )=\lambda
,S_{2}(2,0;\lambda )=0,S_{2}(2,1;\lambda )=\lambda ,\ldots
\end{equation*}%
and%
\begin{equation*}
S_{2}(0,v;\lambda )=\frac{(\lambda -1)^{v}}{v!}.
\end{equation*}

Substituting $\lambda =1$ into (\ref{SN-1}), then we get the Stirling
numbers of the second kind%
\begin{equation*}
S_{2}(n,v)=S_{2}(n,v;1).
\end{equation*}%
(\textit{cf}. \cite{Charamb}-\cite{SrivastavaLiu}; see also the references
cited in each of these earlier works).

\section{Genearlized Apostol-Daehee numbers attached to Dirichlet character $%
\protect\chi $}

In this section, by using the bosonic $p$-adic integral on $\mathbb{Z}_{p}$,
we construct generating functions for the genearlized Apostol-Daehee numbers
attached to Dirichlet character $\chi $. We give relations between these
numbers, the $\lambda $-Bernoulli numbers and the Stirling numbers. We also
give bosonic integral represantation of these numbers. Firstly we give some
standart notations for the Volkenborn integral. Let $\mathbb{Z}_{p}$ be a
set of $p$-adic integers. Let $\ \mathbb{K}$ be a field with a complete
valuation and $C^{1}(\mathbb{Z}_{p}\rightarrow \mathbb{K)}$ be a set of
continuous derivative functions. That is $C^{1}(\mathbb{Z}_{p}\rightarrow 
\mathbb{K)}$ is contained in $\left\{ f:\mathbb{X}\rightarrow \mathbb{K}:f(x)%
\text{ is differentiable and }\frac{d}{dx}f(x)\text{ is continuous}\right\} $%
.

\begin{definition}
(\cite[p. 167-Definition 55.1]{Schikof})The Volkenborn integral of $f\in
C^{1}(\mathbb{Z}_{p}\rightarrow \mathbb{K)}$ is%
\begin{equation}
\int\limits_{\mathbb{Z}_{p}}f\left( x\right) d\mu _{1}\left( x\right) =%
\underset{N\rightarrow \infty }{\lim }\frac{1}{p^{N}}\sum_{x=0}^{p^{N}-1}f%
\left( x\right) .  \label{M}
\end{equation}
\end{definition}

We observe that\ $\mu _{1}\left( x\right) =\mu _{1}\left( x+p^{N}\mathbb{Z}%
_{p}\right) $ is a Haar distrubition, defined by%
\begin{equation*}
\mu _{1}\left( x+p^{N}\mathbb{Z}_{p}\right) =\frac{1}{p^{N}}
\end{equation*}
(\textit{cf}. \cite{Schikof}, \cite{T. Kim}; see also the references cited
in each of these earlier works). In work of Kim \cite{T. Kim}, the
Volkenborn integral is also so-called the\textit{\ bosonic} $p$-adic
Volkenborn integral on $\mathbb{Z}_{p}$. We recall some properties of this
integral as follows:

The Volkenborn integral in terms of the Mahler coefficients is given by the
following formula:%
\begin{equation*}
\int\limits_{\mathbb{Z}_{p}}f\left( x\right) d\mu _{1}\left( x\right)
=\dsum\limits_{n=0}^{\infty }\frac{(-1)^{n}}{n+1}a_{n},
\end{equation*}%
where%
\begin{equation*}
f\left( x\right) =\dsum\limits_{n=0}^{\infty }a_{n}\left( 
\begin{array}{c}
x \\ 
j%
\end{array}%
\right) \in C^{1}(\mathbb{Z}_{p}\rightarrow \mathbb{K)}
\end{equation*}%
(\textit{cf}. \cite[p. 168-Proposition 55.3]{Schikof}). Let $f:\mathbb{Z}%
_{p}\rightarrow \mathbb{K}$ be an analytic function and $f\left( x\right)
=\dsum\limits_{n=0}^{\infty }a_{n}x^{n}$ with $x\in \mathbb{Z}_{p}$. The
Volkenborn integral of this analytic function analytic function is given by%
\begin{equation*}
\int\limits_{\mathbb{Z}_{p}}\left( \dsum\limits_{n=0}^{\infty
}a_{n}x^{n}\right) d\mu _{1}\left( x\right) =\dsum\limits_{n=0}^{\infty
}a_{n}\int\limits_{\mathbb{Z}_{p}}x^{n}d\mu _{1}\left( x\right)
\end{equation*}%
(\textit{cf}. \cite[p. 168-Proposition 55.4]{Schikof}).

The following property is very important to our new results on special
numbers:%
\begin{equation}
\int\limits_{\mathbb{Z}_{p}}f(x+m)d\mu _{1}\left( x\right) =\int\limits_{%
\mathbb{Z}_{p}}f(x)d\mu _{1}\left( x\right) +\dsum\limits_{x=0}^{m-1}\frac{d%
}{dx}f(x)  \label{vi-1}
\end{equation}%
(\textit{cf}. \cite{T. Kim}, \cite{Kim2006TMIC}, \cite{Schikof}, \cite%
{wikipe}; see also the references cited in each of these earlier works).

The $p$-adic $q$-Volkenborn integral was defined by Kim \cite{T. Kim}. Here $%
p$ is a fixed prime. The distribution on $\mathbb{Z}_{p}$ is given by%
\begin{equation*}
\mu _{q}(x+p^{N}\mathbb{Z}_{p})=\frac{q^{x}}{\left[ p^{N}\right] _{q}},
\end{equation*}%
where $q\in \mathbb{C}_{p}$ with $\mid 1-q\mid _{p}<1$ (\textit{cf}. \cite%
{T. Kim}). The $p$-adic $q$-integral of a function $f\in C^{1}(\mathbb{Z}%
_{p}\rightarrow \mathbb{K)}$ is defined by Kim \cite{T. Kim} as follows:%
\begin{equation*}
\int_{\mathbb{Z}_{p}}f(x)d\mu _{q}(x)=\lim_{N\rightarrow \infty }\frac{1}{%
[p^{N}]_{q}}\sum_{x=0}^{p^{N}-1}f(x)q^{x},
\end{equation*}%
where%
\begin{equation*}
\left[ x\right] =\left[ x:q\right] =\left\{ 
\begin{array}{c}
\frac{1-q^{x}}{1-q},q\neq 1 \\ 
x,q=1%
\end{array}%
\right. .
\end{equation*}%
Observe that%
\begin{equation*}
\lim_{q\rightarrow 1}\left[ x:q\right] =x.
\end{equation*}

The Witt's formula for the Bernoulli numbers and polynomials are given as
follows, respectively%
\begin{equation}
\int\limits_{\mathbb{Z}_{p}}x^{n}d\mu _{1}\left( x\right) =B_{n}  \label{M1}
\end{equation}%
and%
\begin{equation}
\int\limits_{\mathbb{Z}_{p}}\left( z+x\right) ^{n}d\mu _{1}\left( x\right)
=B_{n}(z)  \label{wb}
\end{equation}%
(\textit{cf}. \cite{T. Kim}, \cite{Kim2006TMIC}, \cite{Schikof}; see also
the references cited in each of these earlier works).

Let $f\in C^{1}(\mathbb{Z}_{p}\rightarrow \mathbb{K)}$ and 
\begin{equation*}
E^{d}f(x)=f(x+d).
\end{equation*}

Kim \cite[Theorem 1]{KIMaml2008} defined the following functional equation
for the fermionic $p$-adic Volkenborn integral on $\mathbb{Z}_{p}$ as
follows:%
\begin{equation}
q^{n}\int\limits_{\mathbb{Z}_{p}}E^{d}f\left( x\right) d\mu _{q}\left(
x\right) -\int\limits_{\mathbb{Z}_{p}}f\left( x\right) d\mu _{q}\left(
x\right) =\frac{q-1}{\log q}\left( \sum_{j=0}^{n-1}q^{j}f^{^{\prime
}}(j)+\log q\sum_{j=0}^{n-1}q^{j}f(j)\right) ,  \label{TKint.q}
\end{equation}%
where $n$ is a positive integer and%
\begin{equation*}
f^{^{\prime }}(j)=\frac{d}{dx}f(x)\left\vert _{x=j}\right. .
\end{equation*}

\begin{theorem}
\label{TheoremShcrikof}%
\begin{equation}
\int\limits_{\mathbb{Z}_{p}}\left( 
\begin{array}{c}
x \\ 
j%
\end{array}%
\right) d\mu _{1}\left( x\right) =\frac{(-1)^{j}}{j+1}.  \label{C7}
\end{equation}
\end{theorem}

Theorem \ref{TheoremShcrikof} was proved by Schikhof \cite{Schikof}.

Let $\chi $ be a non-trivial Dirichlet character with conductor $d$. Let $%
\lambda $ be a complex number. We set%
\begin{equation}
f(x,t;\lambda )=\lambda ^{x}(1+\lambda t)^{x}\chi (x).  \label{Da-0A}
\end{equation}%
Substituting (\ref{Da-0A}) into (\ref{TKint.q}), we get%
\begin{eqnarray}
\int\limits_{\mathbb{Z}_{p}}\lambda ^{x}(1+\lambda t)^{x}\chi (x)d\mu
_{q}\left( x\right) &=&\frac{q-1}{\log q\left( \left( \lambda q\right)
^{d}(1+\lambda t)^{d}-1\right) }\sum_{j=0}^{d-1}\left( \lambda q\right)
^{j}(1+\lambda t)^{j}\chi (j)\log \left( \lambda +\lambda ^{2}t\right) 
\notag \\
&&+\frac{q-1}{\left( \lambda q\right) ^{d}(1+\lambda t)^{d}-1}%
\sum_{j=0}^{d-1}\left( \lambda q\right) ^{j}(1+\lambda t)^{j}\chi (j).
\label{Da-0}
\end{eqnarray}

From the above integral equation, we define the following generating
function for genearlized Apostol-Daehee numbers attached to Dirichlet
character $\chi $ with conductor $d$ as follows:%
\begin{equation*}
F_{\mathfrak{D}}(t;q,\lambda ,\chi )=\frac{\left( q-1\right) \log \left(
\lambda +\lambda ^{2}t\right) +(q-1)\log q}{\log q}\sum_{j=0}^{d-1}\frac{%
\left( \lambda q\right) ^{j}(1+\lambda t)^{j}\chi (j)}{\left( \lambda
q(1+\lambda t)\right) ^{d}-1}
\end{equation*}%
or%
\begin{equation}
F_{\mathfrak{D}}(t;q,\lambda ,\chi )=\sum_{n=0}^{\infty }\mathfrak{D}%
_{n,\chi }(\lambda ,q)\frac{t^{n}}{n!}.  \label{Da-1}
\end{equation}%
The generalized Apostol-Daehee polynomials attached to Dirichlet character $%
\chi $ are defined by means of the following generating function:(\ref{Da-1})%
\begin{equation}
F_{\mathfrak{D}}(z,t;q,\lambda ,\chi )=F_{\mathfrak{D}}(t;q,\lambda ,\chi
)(1+\lambda t)^{z}=\sum_{n=0}^{\infty }\mathfrak{D}_{n,\chi }(z;\lambda ,q)%
\frac{t^{n}}{n!}.  \label{YY1}
\end{equation}%
We assume that $\left\vert \lambda t\right\vert <1$. Combining the above
function with (\ref{Da-1}), we have%
\begin{equation*}
\sum_{n=0}^{\infty }\mathfrak{D}_{n,\chi }(z;\lambda ,q)\frac{t^{n}}{n!}%
=\sum_{n=0}^{\infty }\left( 
\begin{array}{c}
z \\ 
n%
\end{array}%
\right) \lambda ^{n}\frac{t^{n}}{n!}\sum_{n=0}^{\infty }\mathfrak{D}_{n,\chi
}(\lambda ,q)\frac{t^{n}}{n!}.
\end{equation*}%
Therefore%
\begin{equation*}
\sum_{n=0}^{\infty }\mathfrak{D}_{n,\chi }(z;\lambda ,q)\frac{t^{n}}{n!}%
=\sum_{n=0}^{\infty }\sum_{j=0}^{n}\left( 
\begin{array}{c}
n \\ 
j%
\end{array}%
\right) \lambda ^{n-j}(z)_{n-j}\mathfrak{D}_{j,\chi }(\lambda ,q)\frac{t^{n}%
}{n!}.
\end{equation*}%
Comparing the coefficients of $\frac{t^{m}}{m!}$ on both sides of the above
equation, we arrive at the following theorem:

\begin{theorem}
Let $n\in \mathbb{N}_{0}$. Then we have%
\begin{equation*}
\mathfrak{D}_{n,\chi }(z;\lambda ,q)=\sum_{j=0}^{n}\left( 
\begin{array}{c}
n \\ 
j%
\end{array}%
\right) \lambda ^{n-j}(z)_{n-j}\mathfrak{D}_{j,\chi }(\lambda ,q).
\end{equation*}
\end{theorem}

\begin{remark}
If $q\rightarrow 1$ and $\lambda \rightarrow 1$ and $\chi \equiv 1$, then (%
\ref{YY1}) reduces to 
\begin{equation*}
\frac{\log \left( 1+t\right) }{t}(1+t)^{z}=\sum_{n=0}^{\infty }D_{n}(z)\frac{%
t^{n}}{n!}
\end{equation*}%
(cf. \cite{DSkimDaehee}, \cite{DSKIMopenMath}, \cite{KimITSFdahee}, \cite%
{ysimsek Ascm}, \cite{Simsek2019}).
\end{remark}

We modify (\ref{Da-1}) as follows:%
\begin{eqnarray*}
F_{\mathfrak{D}}(t;q,\lambda ,\chi ) &=&\left( \frac{\left( q-1\right) \log
\lambda }{\log q}+q-1\right) \frac{F_{b2}(\lambda t)}{d\lambda t}%
\sum_{j=0}^{d-1}\left( \lambda q\right) ^{j}\chi (j)F_{A}\left( d\log \left(
1+\lambda t\right) ,\frac{j}{d};\left( \lambda q\right) ^{d}\right) \\
&&+\frac{\left( q-1\right) }{d\log q}\sum_{j=0}^{d-1}\left( \lambda q\right)
^{j}\chi (j)F_{A}\left( d\log \left( 1+\lambda t\right) ,\frac{j}{d};\left(
\lambda q\right) ^{d}\right) .
\end{eqnarray*}%
Combining the above functional equation with (\ref{Ap.B}) and (\ref{Be-1t}),
we get%
\begin{eqnarray*}
\sum_{m=0}^{\infty }m\mathfrak{D}_{m-1,\chi }(\lambda ,q)\frac{t^{m}}{m!}
&=&\left( \frac{\left( q-1\right) \log \lambda }{\log q}+q-1\right)
\sum_{m=0}^{\infty }\sum_{j=0}^{d-1}q^{j}\chi (j)\sum_{l=0}^{m}\left( 
\begin{array}{c}
m \\ 
l%
\end{array}%
\right) \lambda ^{m+j-l-1}b_{m-j}(0) \\
&&\times \sum_{n=0}^{l}d^{n-1}\mathcal{B}_{n}\left( \frac{j}{d};\left(
\lambda q\right) ^{d}\right) S_{1}(l,n)\frac{t^{m}}{m!} \\
&&+\frac{q-1}{\log q}\sum_{m=0}^{\infty }m\sum_{j=0}^{d-1}\left( \lambda
q\right) ^{j}\chi (j)\sum_{n=0}^{m-1}d^{n-1}\mathcal{B}_{n}\left( \frac{j}{d}%
;\left( \lambda q\right) ^{d}\right) S_{1}(m-1,n)\frac{t^{m}}{m!}.
\end{eqnarray*}%
Comparing the coe cients of $\frac{t^{m}}{m!}$ on both sides of the above
equation, we arrive at the following theorem:

\begin{theorem}
Let $m\in \mathbb{N}$. Then we have%
\begin{eqnarray}
\mathfrak{D}_{m-1,\chi }(\lambda ,q) &=&\left( \frac{\left( q-1\right) \log
\lambda }{\log q^{m}}+\frac{q-1}{m}\right) \sum_{j=0}^{d-1}q^{j}\chi
(j)\sum_{l=0}^{m}\left( 
\begin{array}{c}
m \\ 
l%
\end{array}%
\right) \lambda ^{m+j-l-1}b_{m-l}(0)  \notag \\
&&\times \sum_{n=0}^{l}d^{n-1}\mathcal{B}_{n}\left( \frac{j}{d};\left(
\lambda q\right) ^{d}\right) S_{1}(l,n)  \label{Da-2} \\
&&+\frac{q-1}{\log q}\sum_{j=0}^{d-1}\left( \lambda q\right) ^{j}\chi
(j)\sum_{n=0}^{m-1}d^{n-1}\mathcal{B}_{n}\left( \frac{j}{d};\left( \lambda
q\right) ^{d}\right) S_{1}(m-1,n).  \notag
\end{eqnarray}
\end{theorem}

If $q\rightarrow 1$ in (\ref{Da-2}), we get the following corollary:

\begin{corollary}
Let $m\in \mathbb{N}$. Then we have%
\begin{eqnarray}
\mathfrak{D}_{m-1,\chi }(\lambda ) &=&\frac{\log \lambda }{m}%
\sum_{j=0}^{d-1}\chi (j)\sum_{l=0}^{m}\left( 
\begin{array}{c}
m \\ 
l%
\end{array}%
\right) \lambda ^{m+j-l-1}b_{m-l}(0)  \label{Da-3} \\
&&\times \sum_{n=0}^{l}d^{n-1}\mathcal{B}_{n}\left( \frac{j}{d};\lambda
^{d}\right) S_{1}(l,n)+\sum_{j=0}^{d-1}\lambda ^{j}\chi (j)  \notag \\
&&\times \sum_{n=0}^{m-1}d^{n-1}\mathcal{B}_{n}\left( \frac{j}{d};\lambda
^{d}\right) S_{1}(m-1,n).  \notag
\end{eqnarray}
\end{corollary}

If $\lambda =1$, we get the following corollary:

\begin{corollary}
Let $m\in \mathbb{N}_{0}$. Then we have%
\begin{equation}
\mathfrak{D}_{m,\chi }=\sum_{j=0}^{d-1}\chi
(j)\sum_{n=0}^{m}d^{n-1}B_{n}\left( \frac{j}{d}\right) S_{1}(m,n).
\label{Da-4}
\end{equation}
\end{corollary}

By combining (\ref{Da-4}) with the following well-known formula for the
generalized Bernoulli numbers%
\begin{equation*}
B_{n,\chi }=d^{n-1}\sum_{j=0}^{d-1}\chi (j)B_{n}\left( \frac{j}{d}\right) ,
\end{equation*}%
we arrive at the following result.

\begin{corollary}
Let $m\in \mathbb{N}$. Then we have%
\begin{equation}
\mathfrak{D}_{m,\chi }=\sum_{n=0}^{m}B_{n,\chi }S_{1}(m,n).  \label{Da-5}
\end{equation}
\end{corollary}

Assume that $\left\vert \lambda t\right\vert <1$. Combining (\ref{Da-0})
with (\ref{Da-1}), we get%
\begin{equation*}
\sum_{n=0}^{\infty }\left( \int\limits_{\mathbb{Z}_{p}}(x)_{n}\lambda
^{x+n}\chi (x)d\mu _{q}\left( x\right) \right) \frac{t^{n}}{n!}%
=\sum_{n=0}^{\infty }\mathfrak{D}_{n,\chi }(\lambda ,q)\frac{t^{n}}{n!}
\end{equation*}%
Comparing the coefficients of $\frac{t^{n}}{n!}$ on both sides of the above
equation, we arrive at the following theorem:

\begin{theorem}
Let $n\in \mathbb{N}_{0}$. Then we have%
\begin{equation}
\mathfrak{D}_{n,\chi }(\lambda ,q)=\int\limits_{\mathbb{Z}%
_{p}}(x)_{n}\lambda ^{x+n}\chi (x)d\mu _{q}\left( x\right) .  \label{Da-i1}
\end{equation}
\end{theorem}

We define genearlized Apostol-Daehee polynomials attached to Dirichlet
character $\chi $ with conductor $d$ as follows:%
\begin{equation*}
G(t,x;q,\lambda ,\chi )=(1+\lambda t)^{x}F_{\mathfrak{D}}(t;q,\lambda ,\chi
)=\sum_{n=0}^{\infty }\mathfrak{D}_{n,\chi }(x;\lambda ,q)\frac{t^{n}}{n!}.
\end{equation*}%
From this equation, we get the following theorem:

\begin{theorem}
Let $n\in \mathbb{N}_{0}$. Then we have%
\begin{equation*}
\mathfrak{D}_{n,\chi }(x;\lambda ,q)=\sum_{j=0}^{n}\left( 
\begin{array}{c}
n \\ 
j%
\end{array}%
\right) \lambda ^{j}(x)_{j}\mathfrak{D}_{n-j,\chi }(\lambda ,q).
\end{equation*}
\end{theorem}

Substituting $\lambda t=e^{t}-1$ into (\ref{Da-1}), we get%
\begin{equation*}
\frac{\left( q-1\right) \left( \log \lambda +t\right) +\left( q-1\right)
\log q}{\log q}\sum_{j=0}^{d-1}\frac{\left( \lambda q\right) ^{j}e^{tj}\chi
(j)}{\left( \lambda q\right) ^{d}e^{td}-1}=\sum_{n=0}^{\infty }\frac{%
\mathfrak{D}_{n,\chi }(\lambda ,q)}{\lambda ^{n}}\frac{\left( e^{t}-1\right)
^{n}}{n!}.
\end{equation*}%
By substituting (\ref{SN-1}) into the above equation, 
\begin{equation}
\frac{\left( q-1\right) \left( \log \lambda +\log q+t\right) }{\log q}%
\sum_{j=0}^{d-1}\frac{\left( \lambda q\right) ^{j}e^{tj}\chi (j)}{\left(
\lambda q\right) ^{d}e^{td}-1}=\sum_{n=0}^{\infty }\frac{\mathfrak{D}%
_{n,\chi }(\lambda ,q)}{\lambda ^{n}}\sum_{m=0}^{\infty }S_{2}(m,n)\frac{%
t^{m}}{m!}.  \label{FEN-1}
\end{equation}%
Since $S_{2}(m,n)=0$ if $n>m$, we get%
\begin{eqnarray*}
&&\sum_{m=0}^{\infty }m\sum_{n=0}^{m-1}\frac{\mathfrak{D}_{n,\chi }(\lambda
,q)S_{2}(m-1,n)}{\lambda ^{n}}\frac{t^{m}}{m!} \\
&=&\sum_{m=0}^{\infty }md^{m-1}\frac{q-1}{\log q}\log \left( \lambda
q\right) \sum_{j=0}^{d-1}\left( \lambda q\right) ^{j}\chi (j)\mathcal{B}%
_{n}\left( \frac{j}{d};\left( \lambda q\right) ^{d}\right) \frac{t^{m}}{m!}
\\
&&+\sum_{m=0}^{\infty }d^{m-1}\frac{q-1}{m\log q}\sum_{j=0}^{d-1}\left(
\lambda q\right) ^{j}\chi (j)\mathcal{B}_{n}\left( \frac{j}{d};\left(
\lambda q\right) ^{d}\right) \frac{t^{m}}{m!}.
\end{eqnarray*}%
Comparing the coefficients of $\frac{t^{m}}{m!}$ on both sides of the above
equation, we arrive at the following theorem:

\begin{theorem}
Let Let $m\in \mathbb{N}$. Then we have%
\begin{eqnarray}
\sum_{n=0}^{m-1}\frac{\mathfrak{D}_{n,\chi }(\lambda ,q)S_{2}(m-1,n)}{%
\lambda ^{n}} &=&\frac{d^{m-1}\left( q-1\right) \log \left( \lambda q\right) 
}{\log q}\sum_{j=0}^{d-1}\left( \lambda q\right) ^{j}\chi (j)\mathcal{B}%
_{m-1}\left( \frac{j}{d};\left( \lambda q\right) ^{d}\right)  \notag \\
&&+\frac{\left( q-1\right) d^{m-1}}{m\log q}\sum_{j=0}^{d-1}\left( \lambda
q\right) ^{j}\chi (j)\mathcal{B}_{m}\left( \frac{j}{d};\left( \lambda
q\right) ^{d}\right) .  \label{Da-6}
\end{eqnarray}
\end{theorem}

If $\lambda =1$ and $q\rightarrow 1$ in (\ref{Da-6}), we get the following
corollary:

\begin{corollary}
Let $m\in \mathbb{N}$. Then we have%
\begin{equation*}
\sum_{n=0}^{m}\mathfrak{D}_{n,\chi }S_{2}(m,n)=d^{m-1}\sum_{j=0}^{d-1}\chi
(j)B_{m}\left( \frac{j}{d}\right)
\end{equation*}%
or%
\begin{equation*}
B_{m,\chi }=\sum_{n=0}^{m}\mathfrak{D}_{n,\chi }S_{2}(m,n).
\end{equation*}
\end{corollary}

By using (\ref{FEN-1}) with $S_{2}(m,n)=0$ if $n>m$, we get%
\begin{equation*}
\frac{\left( q-1\right) \left( \log \lambda +\log q+t\right) }{\left(
\lambda ^{d}q^{d}-1\right) \log q}\sum_{j=0}^{d-1}\frac{\left( \lambda
q\right) ^{j}e^{tj}\chi (j)\left( 1-\frac{1}{\lambda ^{d}q^{d}}\right) }{%
e^{td}-\frac{1}{\lambda ^{d}q^{d}}}=\sum_{m=0}^{\infty }\sum_{n=0}^{m}\frac{%
\mathfrak{D}_{n,\chi }(\lambda ,q)S_{2}(m,n)}{\lambda ^{n}}\frac{t^{m}}{m!}.
\end{equation*}%
By substituting (\ref{FEN}) into the above equation, since $S_{2}(m,n)=0$ if 
$n>m$, we get%
\begin{eqnarray*}
&&\sum_{m=0}^{\infty }\sum_{n=0}^{m}\frac{\mathfrak{D}_{n,\chi }(\lambda
,q)S_{2}(m,n)}{\lambda ^{n}}\frac{t^{m}}{m!} \\
&=&\frac{\left( q-1\right) \left( \log \lambda +\log q\right) }{\left(
\lambda ^{d}q^{d}-1\right) \log q}\sum_{m=0}^{\infty
}d^{m}\sum_{j=0}^{d-1}\left( \lambda q\right) ^{j}\chi (j)H_{m}\left( \frac{j%
}{d};\frac{1}{\lambda ^{d}q^{d}}\right) \frac{t^{m}}{m!} \\
&&+\sum_{m=0}^{\infty }md^{m-1}\frac{q-1}{\log q}\sum_{j=0}^{d-1}\left(
\lambda q\right) ^{j}\chi (j)H_{m-1}\left( \frac{j}{d};\frac{1}{\lambda
^{d}q^{d}}\right) \frac{t^{m}}{m!}.
\end{eqnarray*}%
Comparing the coefficients of $\frac{t^{m}}{m!}$ on both sides of the above
equation, we arrive at the following theorem:

\begin{theorem}
\begin{eqnarray}
\sum_{n=0}^{m}\frac{\mathfrak{D}_{n,\chi }(\lambda ,q)S_{2}(m,n)}{\lambda
^{n}} &=&\frac{d^{m}\left( q-1\right) \left( \log \lambda +\log q\right) }{%
\left( \lambda ^{d}q^{d}-1\right) \log q}\sum_{j=0}^{d-1}\left( \lambda
q\right) ^{j}\chi (j)H_{m}\left( \frac{j}{d};\frac{1}{\lambda ^{d}q^{d}}%
\right)  \notag \\
&&+md^{m-1}\frac{q-1}{\log q}\sum_{j=0}^{d-1}\left( \lambda q\right)
^{j}\chi (j)H_{m-1}\left( \frac{j}{d};\frac{1}{\lambda ^{d}q^{d}}\right) .
\label{FEN2}
\end{eqnarray}
\end{theorem}

If $q\rightarrow 1$ in (\ref{FEN2}), we get the following corollary:

\begin{corollary}
\begin{equation*}
\sum_{n=0}^{m}\frac{\mathfrak{D}_{n,\chi }(\lambda ,1)S_{2}(m,n)}{\lambda
^{n}}=md^{m-1}\sum_{j=0}^{d-1}\lambda ^{j}\chi (j)H_{m-1}\left( \frac{j}{d};%
\frac{1}{\lambda ^{d}}\right) .
\end{equation*}
\end{corollary}

\section{Genearlized Apostol-Changhee numbers attached to Dirichlet
character $\protect\chi $}

In this section, by using the\textit{\ }fermionic $p$-adic integral on $%
\mathbb{Z}_{p}$, we construct generating functions for the genearlized
Apostol-Changhee numbers and polynomials attached to Dirichlet character $%
\chi $. We give relations between these numbers, the Apostol-Euler numbers
and the Stirling numbers. We also give bosonic integral represantation of
these numbers. Firstly we give some standart notations for the fermionic $p$%
-adic integral on $\mathbb{Z}_{p}$.

The \textit{fermionic} $p$-adic Volkenborn integral on $\mathbb{Z}_{p}$ is
given by 
\begin{equation}
\int\limits_{\mathbb{Z}_{p}}f\left( x\right) d\mu _{-1}\left( x\right) =%
\underset{N\rightarrow \infty }{\lim }\sum_{x=0}^{p^{N}-1}\left( -1\right)
^{x}f\left( x\right)  \label{Mmm}
\end{equation}%
where%
\begin{equation*}
\mu _{1}\left( x+p^{N}\mathbb{Z}_{p}\right) =\frac{(-1)^{x}}{p^{N}}
\end{equation*}%
(\textit{cf}. \cite{Kim2006TMIC}).

Let $f\in C^{1}(\mathbb{Z}_{p}\rightarrow \mathbb{K)}$ and 
\begin{equation*}
E^{d}f(x)=f(x+d).
\end{equation*}%
Kim \cite{KIMjmaa2017} defined the following functional equation for the
fermionic $p$-adic Volkenborn integral on $\mathbb{Z}_{p}$ as follows:

\begin{equation}
\int\limits_{\mathbb{Z}_{p}}E^{d}f\left( x\right) d\mu _{-1}\left( x\right)
-(-1)^{d}\int\limits_{\mathbb{Z}_{p}}f\left( x\right) d\mu _{-1}\left(
x\right) =2\sum_{j=0}^{d-1}(-1)^{d-1-j}f(j),  \label{TKint}
\end{equation}%
where $d$ is a positive integer. Substituting $d=1$ into (\ref{TKint}), we
have%
\begin{equation*}
\int\limits_{\mathbb{Z}_{p}}f\left( x+1\right) d\mu _{-1}\left( x\right)
+\int\limits_{\mathbb{Z}_{p}}f\left( x\right) d\mu _{-1}\left( x\right)
=2f(0)
\end{equation*}%
(\textit{cf}. \cite{KIMjmaa2017}). By using (\ref{Mmm}), the Witt's formula
for the Euler numbers and polynomials are given as follows, respectively%
\begin{equation}
\int\limits_{\mathbb{Z}_{p}}x^{n}d\mu _{-1}\left( x\right) =E_{n}
\label{Mm1}
\end{equation}%
and%
\begin{equation}
\int\limits_{\mathbb{Z}_{p}}\left( z+x\right) ^{n}d\mu _{-1}\left( x\right)
=E_{n}(z)  \label{we}
\end{equation}%
(\textit{cf}. \cite{Kim2006TMIC}, \cite{KIMjang}; see also the references
cited in each of these earlier works).

\begin{theorem}
\label{ThoremKIM}%
\begin{equation}
\int\limits_{\mathbb{Z}_{p}}\left( 
\begin{array}{c}
x \\ 
j%
\end{array}%
\right) d\mu _{-1}\left( x\right) =\frac{(-1)^{j}}{2^{j}}.  \label{est-3}
\end{equation}
\end{theorem}

Theorem \ref{ThoremKIM} was proved by Kim \textit{et al.} \cite{DSkim2}.

In \cite{KIMaml2008}, Kim gave $q$-analoges of (\ref{TKint}) as follows:%
\begin{equation}
q^{d}\int\limits_{\mathbb{Z}_{p}}E^{d}f\left( x\right) d\mu _{-q}\left(
x\right) +\int\limits_{\mathbb{Z}_{p}}f\left( x\right) d\mu _{-q}\left(
x\right) =\left[ 2\right] \sum_{j=0}^{d-1}(-1)^{l}q^{j}f(j),
\label{TkimFint}
\end{equation}%
where $d$ is an positif odd integer.

Substituting (\ref{Da-0A}) into (\ref{TkimFint}), we get%
\begin{equation}
\int\limits_{\mathbb{Z}_{p}}\lambda ^{x}(1+\lambda t)^{x}\chi (x)d\mu
_{-q}\left( x\right) =\frac{\left[ 2\right] }{\left( \lambda q\right)
^{d}(1+\lambda t)^{d}+1}\sum_{j=0}^{d-1}(-1)^{j}\chi (j)\left( \lambda
q\right) ^{j}(1+\lambda t)^{j}.  \label{F1}
\end{equation}%
By using this equation, we define generating functions for the generalized
Apostol-Changhee type numbers and polynomials by means of the following
generating functions, respectively:%
\begin{equation}
F_{\mathfrak{E}}(t;\lambda ,q,\chi )=\frac{\left[ 2\right] }{\left( \lambda
q\right) ^{d}(1+\lambda t)^{d}+1}\sum_{j=0}^{d-1}(-1)^{j}\chi (j)\left(
\lambda q\right) ^{j}(1+\lambda t)^{j},  \label{F2}
\end{equation}

\begin{equation*}
F_{\mathfrak{C}}(t;\lambda ,q,\chi )=\sum_{n=0}^{\infty }\mathfrak{Ch}%
_{n,\chi }(\lambda ,q)\frac{t^{n}}{n!}
\end{equation*}%
and%
\begin{eqnarray}
F_{\mathfrak{E}}(t,z;\lambda ,q,\chi ) &=&F_{\mathfrak{E}}(t;\lambda ,q,\chi
)(1+\lambda t)^{z}  \label{YY-4} \\
&=&\sum_{n=0}^{\infty }\mathfrak{Ch}_{n,\chi }(z;\lambda ,q)\frac{t^{n}}{n!}.
\notag
\end{eqnarray}%
We assume that $\left\vert \lambda t\right\vert <1$. Combining (\ref{F2})
with the above equation, we get%
\begin{equation*}
\sum_{n=0}^{\infty }\mathfrak{Ch}_{n,\chi }(z;\lambda ,q)\frac{t^{n}}{n!}%
=\sum_{n=0}^{\infty }\left( 
\begin{array}{c}
z \\ 
n%
\end{array}%
\right) \lambda ^{n}\frac{t^{n}}{n!}\sum_{n=0}^{\infty }\mathfrak{Ch}%
_{n,\chi }(\lambda ,q)\frac{t^{n}}{n!}.
\end{equation*}%
Therefore%
\begin{equation*}
\sum_{n=0}^{\infty }\mathfrak{Ch}_{n,\chi }(z;\lambda ,q)\frac{t^{n}}{n!}%
=\sum_{j=0}^{n}\left( 
\begin{array}{c}
n \\ 
j%
\end{array}%
\right) \lambda ^{n-j}(z)_{n-j}\mathfrak{Ch}_{j,\chi }(\lambda ,q)\frac{t^{n}%
}{n!}.
\end{equation*}%
Comparing the coefficients of $\frac{t^{m}}{m!}$ on both sides of the above
equation, we arrive at the following theorem:

\begin{theorem}
Let $n\in \mathbb{N}_{0}$. Then we have%
\begin{equation*}
\mathfrak{Ch}_{n,\chi }(z;\lambda ,q)=\sum_{j=0}^{n}\left( 
\begin{array}{c}
n \\ 
j%
\end{array}%
\right) \lambda ^{n-j}(z)_{n-j}\mathfrak{Ch}_{j,\chi }(\lambda ,q).
\end{equation*}
\end{theorem}

\begin{remark}
If $q\rightarrow 1$ and $\lambda \rightarrow 1$ and $\chi \equiv 1$, then (%
\ref{YY-4}) reduces to%
\begin{equation*}
\frac{2}{t+2}(1+t)^{z}=\sum_{n=0}^{\infty }Ch_{n}(z)\frac{t^{n}}{n!}
\end{equation*}%
(cf. \cite{DSkim2}, \cite{KimCHAReul}, \cite{Simsek2019}).
\end{remark}

By using (\ref{F1}), we get the following result:

\begin{theorem}
Let $n\in \mathbb{N}_{0}$. Then we have%
\begin{equation}
\mathfrak{Ch}_{n,\chi }(\lambda ,q)=\int\limits_{\mathbb{Z}_{p}}\lambda
^{x+n}(x)_{n}\chi (x)d\mu _{-q}\left( x\right) .  \label{F3b}
\end{equation}
\end{theorem}

By using (\ref{F1}), we derive the following functional equation: 
\begin{equation*}
F_{\mathfrak{C}}(t,x;\lambda ,q,\chi )=\frac{\left[ 2\right] }{2}%
\sum_{j=0}^{d-1}(-1)^{j}\chi (j)\left( \lambda q\right) ^{j}F_{P1}\left(
d\log \left( 1+\lambda t\right) ,\frac{j}{d};1,\left( \lambda q\right)
^{d}\right)
\end{equation*}%
Combining this equation with (\ref{Cad3}) and (\ref{F2}), we get%
\begin{eqnarray*}
&&\frac{\left[ 2\right] }{2}\sum_{j=0}^{d-1}(-1)^{j}\chi (j)\left( \lambda
q\right) ^{j}\sum_{n=0}^{\infty }d^{n}\mathcal{E}_{n}\left( \frac{j}{d}%
;\left( \lambda q\right) ^{d}\right) \frac{\left( \log \left( 1+\lambda
t\right) \right) ^{n}}{n!} \\
&=&\sum_{m=0}^{\infty }\mathfrak{C}_{m,\chi }(\lambda ,q)\frac{t^{m}}{m!}.
\end{eqnarray*}%
Combining (\ref{S1}), we have%
\begin{eqnarray*}
&&\frac{\left[ 2\right] }{2}\sum_{j=0}^{d-1}\chi (j)\left( -\lambda q\right)
^{j}\sum_{n=0}^{\infty }d^{n}\mathcal{E}_{n}\left( \frac{j}{d};\left(
\lambda q\right) ^{d}\right) \sum_{m=0}^{\infty }S_{1}(m,n)\frac{\left(
\lambda t\right) ^{m}}{m!} \\
&=&\sum_{m=0}^{\infty }\mathfrak{C}_{n,\chi }(\lambda ,q)\frac{t^{m}}{m!}.
\end{eqnarray*}%
Comparing the coefficients of $\frac{t^{m}}{m!}$ on both sides of the above
equation, since $S_{1}(m,n)=0$, $m<n$, we arrive at the following theorem:

\begin{theorem}
Let $m$ be a nonnegative integer. Then we have%
\begin{equation}
\mathfrak{Ch}_{m,\chi }(\lambda ,q)=\sum_{j=0}^{d-1}(-q)^{j}\chi
(j)\sum_{n=0}^{m}\lambda ^{j+n}d^{n}\mathcal{E}_{n}\left( \frac{j}{d};\left(
\lambda q\right) ^{d}\right) S_{1}(m,n).  \label{F3-a}
\end{equation}
\end{theorem}

By combining (\ref{F3-a}) with the following well-known formula for the
generalized Euler numbers%
\begin{equation*}
\mathfrak{Ch}_{n,\chi }(\lambda )=d^{n}\sum_{j=0}^{d-1}(-1)^{j}\lambda
^{j}\chi (j)\mathcal{E}_{n}\left( \frac{j}{d};\lambda ^{d}\right) ,
\end{equation*}%
we arrive at the following result.

\begin{corollary}
Let $m\in \mathbb{N}_{0}$. Then we have%
\begin{equation}
\mathfrak{Ch}_{m,\chi }(\lambda ,q)=\sum_{n=0}^{m}\mathcal{E}_{n,\chi
}(q\lambda )S_{1}(m,n).  \label{F3}
\end{equation}
\end{corollary}

Substituting%
\begin{equation*}
\lambda t=e^{t}-1
\end{equation*}%
into (\ref{F2}), we get%
\begin{equation}
\frac{\left[ 2\right] }{\left( \lambda q\right) ^{d}e^{dt}+1}%
\sum_{j=0}^{d-1}(-1)^{j}\chi (j)\left( \lambda q\right)
^{j}e^{jt}=\sum_{n=0}^{\infty }\mathfrak{Ch}_{n,\chi }(\lambda ,q)\frac{%
\left( e^{t}-1\right) ^{n}}{n!}.  \label{FEN-3}
\end{equation}

By substituting (\ref{FEN}) into the above equation, since $S_{2}(m,n)=0$ if 
$n>m$, we get%
\begin{equation*}
\sum_{m=0}^{\infty }\left( \frac{\left[ 2\right] }{2}%
\sum_{j=0}^{d-1}(-1)^{j}\chi (j)\left( \lambda q\right) ^{j}d^{m}\mathcal{E}%
_{m}\left( \frac{j}{d};\left( \lambda q\right) ^{d}\right) \right) \frac{%
t^{m}}{m!}=\sum_{m=0}^{\infty }\sum_{n=0}^{m}\frac{\mathfrak{Ch}_{n,\chi
}(\lambda ,q)S_{2}(m,n)}{\lambda ^{n}}\frac{t^{m}}{m!}.
\end{equation*}%
Comparing the coefficients of $\frac{t^{m}}{m!}$ on both sides of the above
equation, we arrive at the following theorem:

\begin{theorem}
\begin{equation*}
\frac{\left[ 2\right] }{2}\sum_{j=0}^{d-1}(-1)^{j}\chi (j)\left( \lambda
q\right) ^{j}d^{m}\mathcal{E}_{m}\left( \frac{j}{d};\left( \lambda q\right)
^{d}\right) =\sum_{n=0}^{m}\frac{\mathfrak{Ch}_{n,\chi }(\lambda
,q)S_{2}(m,n)}{\lambda ^{n}}
\end{equation*}%
or%
\begin{equation*}
\mathcal{E}_{m},\chi =\frac{2}{\left[ 2\right] }\sum_{n=0}^{m}\frac{%
\mathfrak{Ch}_{n,\chi }(\lambda ,q)S_{2}(m,n)}{\lambda ^{n}}.
\end{equation*}
\end{theorem}

By combining (\ref{FEN}) with (\ref{FEN-3}), we get%
\begin{eqnarray*}
&&\sum_{m=0}^{\infty }\left( \frac{\left[ 2\right] }{\left( \lambda
^{d}q^{d}+1\right) }\sum_{j=0}^{d-1}(-1)^{j}\chi (j)\left( \lambda q\right)
^{j}d^{m}H_{m}\left( \frac{j}{d};-\frac{1}{\left( \lambda q\right) ^{d}}%
\right) \right) \frac{t^{m}}{m!} \\
&=&\sum_{m=0}^{\infty }\sum_{n=0}^{m}\frac{\mathfrak{Ch}_{n,\chi }(\lambda
,q)S_{2}(m,n)}{\lambda ^{n}}\frac{t^{m}}{m!}.
\end{eqnarray*}%
Comparing the coefficients of $\frac{t^{m}}{m!}$ on both sides of the above
equation, we arrive at the following theorem:

\begin{theorem}
\begin{equation*}
\frac{\left[ 2\right] }{\left( \lambda ^{d}q^{d}+1\right) }%
\sum_{j=0}^{d-1}(-1)^{j}\chi (j)\left( \lambda q\right) ^{j}d^{m}H_{m}\left( 
\frac{j}{d};-\frac{1}{\left( \lambda q\right) ^{d}}\right) =\sum_{n=0}^{m}%
\frac{\mathfrak{Ch}_{n,\chi }(\lambda ,q)S_{2}(m,n)}{\lambda ^{n}}.
\end{equation*}
\end{theorem}

\begin{theorem}
Let $d$ be an odd integer. Let $m$ be a positive integer. Then we have%
\begin{equation*}
\mathcal{B}_{m,\chi }(\lambda )=m\sum_{n=0}^{m-1}(-1)^{n}\frac{\mathfrak{Ch}%
_{n,\chi }(-\lambda ,q)S_{2}(m,n)}{\lambda ^{n}}.
\end{equation*}
\end{theorem}

\begin{acknowledgement}
The first author was supported by the \textit{Scientific Research Project
Administration of Akdeniz University.}
\end{acknowledgement}

\end{document}